\documentclass[11pt]{article}
\usepackage{amsmath,amsthm}

\let\workingver=n
\def\refeq#1{\if\workingver y(\ref{#1})-[[#1]]\else(\ref{#1})\fi}
\def\refth#1{\if\workingver y\ref{#1}-[[#1]]\else\ref{#1}\fi}
\def\mylabel#1{\if\workingver y\label{#1}{\bf\ \ [[#1]]\ \ }
\else\label{#1}\fi}
\def\mybibitem#1{\if\workingver y\bibitem{#1}{\bf\ \ [[#1]]\ \ }
\else\bibitem{#1}\fi}

\newtheorem{thm}{Theorem}
\newtheorem{lem}[thm]{Lemma}
\newtheorem{defn}[thm]{Definition}

\newtheorem{rem}[thm]{Remark}
\newtheorem{conj}[thm]{Conjecture}

\newcommand{\bb}{\bar}

\newcommand{\beq}{\begin{equation}}
\newcommand{\gata}{\end{equation}}

\font\tenrm=cmr10

\font\ninerm=cmr9

\def\sms{\small\scshape}

\def\ZZ{{\bf Z}_2^n}

\def\op{\oplus}

\begin{document}
\title{ Fast Evaluation, Weights and Nonlinearity of
Rotation-Symmetric Functions}
\author{Thomas W. Cusick\footnote{
State University of New York at Buffalo, Department of Mathematics,
Buffalo, NY 14260-2900, e-mail: cusick@math.buffalo.edu},\quad
Pantelimon St\u anic\u a\footnote{
Auburn University Montgomery, Department of Mathematics,
Montgomery, AL 36124-4023, e-mail: stanpan@strudel.aum.edu}
\thanks{The second author is on leave from the Institute of
Mathematics of Romanian Academy, Bucharest, Romania}}

\date{August 10, 2000}
\maketitle

\pagestyle{myheadings}

\hrule

\vspace{.3cm}
{\noindent\bf Abstract}

\vspace{.2cm}
{\tenrm
We study the nonlinearity and the weight of the {\em rotation-symmetric}
({\em RotS}) functions defined by Pieprzyk and Qu \cite{PQ}. We give
exact results for the nonlinearity and weight  of 2-degree {\em RotS}
functions with the help of the semi-bent functions \cite{CLK}
 and we give the generating function for the weight of the 3-degree
 {\em RotS} function. Based on the numerical examples and our
 observations we state a conjecture on the nonlinearity and weight of
the 3-degree {\em RotS} function.
}

\vspace{.4cm}

{\ninerm
}

\vspace{.3cm}
\hrule
\pagestyle{myheadings}

\vspace{1.4cm}

\section{Motivation}

\baselineskip=2\baselineskip

Hash functions are used to map a large collection of {\em
messages} into a small set of {\em message digests} and can be
used to generate efficiently both signatures and message
authentication codes, and they can be also used as one-way functions in
key agreement and key establishment protocols.
There are two approaches to the study of { hash functions}:
{\em Information Theory} and {\em Complexity Theory}. The first
method provides unconditional security -- an enemy cannot attack
such systems even if he/she has unlimited computing power.
Unfortunately, this is still a theoretical approach and is
generally impractical \cite{BSP}. In the second method based on
complexity theory, some assumptions are made on the computing
power of the enemy or the weaknesses of the existing systems and
algorithms. The best we can hope for is to estimate the computing
power necessary for the attacker to break the algorithm.
Recent progress in interpolation cryptanalysis \cite{JK} and high
order differential cryptanalysis \cite{SMK} has shown that the
algebraic degree is an important factor in the design of
cryptographic primitives.
In fact, in \cite{SMK} the algebraic degree is the crucial
parameter in determining how secure certain cryptosystems are against
higher order differential attacks.
Together with
propagation, differential and nonlinearity profile, resiliency,
correlation-immunity, local and global avalanche characteristics
they form a class of design criteria which we have to consider in
the design of such primitives.

In \cite{PQ}, Pieprzyk and Qu studied some functions, which they
called {\em ro\-ta\-tion-sym\-me\-tric} ({\em RotS}) as components
in the rounds of a
hashing algorithm. It turns out that the degree-two {\em RotS} function
takes $\frac{3n-1}{2}+6(m-1)$ operations (additions and multiplications)
to evaluate in $m$
consecutive rounds of a hashing algorithm. In
\cite{SMK} the authors showed how to break in less than 20
milli-seconds a block cipher that employs low algebraic degree
(quadratic) Boolean functions as its S-boxes and is provably
secure against linear and differential attacks. Therefore, it is
necessary to employ high degree {\em RotS} functions in our algorithms.
To protect from differential attack, we
need {\em RotS} functions with high nonlinearity.
In this paper we aim to complete the study begun by Pieprzyk and
 Qu \cite{PQ} on
the two-degree {\em RotS} functions and we construct
the three-degree {\em RotS} functions and  we prove some results
about their weights
and nonlinearity.

\section{Preliminaries}

Let $n\geq 6$ be a positive integer and $W_n=\{0,1\}^n$ be the space of
binary vectors.
Denote $\alpha_0=(0,\ldots, 0,0),\alpha_1=(0,\ldots,0,1),
\ldots,\alpha_{2^n-1}=(1,\ldots,1).$
We use the lexicographical order on the sequence $\alpha$, that is
$\alpha_0<\alpha_1<\cdots<\alpha_{2^n-1}.$ The Boolean functions
will be written in their algebraic normal form (when
$\alpha=(a_1,\ldots,a_n))$ as
\[
f(x)=\oplus_{\alpha\in W_n} c_{\alpha} x_1^{a_1}\cdots x_n^{a_n},
\]
where $c_{\alpha}\in W_1$. The truth table of $f$ is the binary
sequence
\beq
\label{f-blocks}
f=(v_1,v_2,\ldots,v_{2^n}),
\end{equation}
where the bits $v_1=f\bigl ( (0,\ldots,0)\bigr )$,
$v_2=f\bigl ( (0,\ldots,0,1)\bigr ),\ldots$.
We shall identify the function $f$ with its
vector representation in (\ref{f-blocks}).
We call a function {\em balanced} if the number of ones is equal
to the number of zeroes in its truth table. The {\em Hamming
weight} of a binary vector $v$, denoted by $wt(v)$ is defined as the
number of ones it contains. The {\em Hamming distance} between two
functions $f,g:W_n\to W_1$, denoted by $d(f,g)$ is defined as
$wt(f\oplus g)$.
The nonlinearity of a function $f$, denoted by $N_f$ is defined as
\[
\min_{\phi\in A_n} d(f,\phi),
\]
where $A_n$ is the class of all affine function on $W_n$.
We say that $f$ satisfies the {\em propagation criterion (PC)} with
respect to $c$ if
\beq
\label{sac}
\sum_{x\in \ZZ} f(x)\oplus f(x\op c)=2^{n-1}.
\end{equation}
If $f$ satisfies the PC with respect to all vectors of weight 1, $f$
is called an {\em SAC} ({\em Strict Avalanche Criterion}) function.
If the above relation happens for any $c$ with $wt(c)\leq s$, we
say that $f$ satisfies
$PC(s)$, and if $s=n$, then we say that $f$ is a {\em bent function}.
If two functions $g,h$, on $W_n$, satisfy
$g(x)=h(Ax\op a)\op (b\cdot x)\op c$
with
$a,b\in W_n, c\in W_1$, and $A$ a $2k\times 2k$
nonsingular matrix, we say that $g$ is {\em affinely equivalent} to $h$.

\begin{defn}
The class of {\em rotation-symmetric} ({\em RotS}) functions includes all
Boolean functions $f:W_n\to W_1$ such that $f(x_1,\ldots,x_n)=
f(\rho(x_1),\ldots,\\
\rho(x_n)),$
where $\rho(x_i)=x_{i+1}$, and $x_{n+1}:=x_1$.
\end{defn}

\def\hF{{\hat {\cal F}}}

As in \cite{PQ}, we denote by $\rho$ the permutation $\rho(i)=i+1,
\rho(n)=1$. By abuse of notation we use the same letter for the
transformation which acts on each variable by $\rho(x_i)=x_{i+1},
\rho(x_n)=x_1.$
By $\hat g$ we mean $(-1)^g$.
We define the {\em Walsh-Hadamard transform} of a $g\in W_n$ to be the map
$\hat {\cal F}_{\hat g}:W_n\to {\bf R}$,
\[
\hF_{\hat g}(w)=\sum_{x\in W_n} \hat g(x) (-1)^{w\cdot x}.
\]
The {\em correlation value} between $g$ and $h$ it is defined by
\[
c(g,h)=1-\frac{d(g,h)}{2^{n-1}}.
\]

If $U$ is a string of bits, then $\bar U$ denotes the complemented
string with 0 and 1 interchanged. If $X$ is a 4-bit
block or a string of blocks,
by $(X)_u$ or $X_u$ we shall mean the string obtained by
concatenation of $u$
copies of $X$. The concatenation of two strings $u,v$ will be
denoted by $u v$ or $u||v$.
 Now we define two sets of 4-bit strings
\begin{eqnarray*}
T_1=\{A&=0,0,1,1;\ {\bb A}=1,1,0,0;\ B=0,1,0,1;{\bb B}=1,0,1,0;\\
C&=0,1,1,0;\ {\bb C}=1,0,0,1;\ D=0,0,0,0;{\bb D}=1,1,1,1\}
\end{eqnarray*}
and
\begin{eqnarray*}
T_2=\{
U&=1,0,0,0;\ {\bb U}=0,1,1,1;\ V=0,0,0,1;\ {\bb V}=1,1,1,0;\\
X&=0,1,0,0;\ {\bb X}=1,0,1,1;\ Y=0,0,1,0;\ {\bb Y}=1,1,0,1 \}.
\end{eqnarray*}


\section{The second degree rotation-symmetric function}


In \cite{PQ} the authors proved that the homogeneous rotation symmetric
function of degree 2,  $f_2=x_1x_l+x_2x_{l+1}+\cdots x_n x_{n+l-1},$
(the subscript $w$ is taken as $((w-1)\bmod n)+1$) has  good
nonlinearity and good avalanche properties. Precisely, they proved
\begin{thm}
The function $f_2$ has the following properties:
\begin{itemize}
\item[(i)] the Hamming weight satisfies
$2^{n-2}\leq wt(f)\leq 2^n-2^{n-2}$,
\item[(ii)] the nonlinearity satisfies $N_{f_2}\geq 2^{n-2}$,
\item[(iii)] if $n$ is odd, then $N_{f_2}=2^{n-1}-2^{\frac{n-1}{2}}$
and $f_2$
is balanced,
\item[(iv)] the function satisfies the PC with respect to all
 vectors $\alpha$ of
weight $0<wt(\alpha)<n.$ In particular $f$ is an SAC function.
\end{itemize}
\end{thm}

In the same paper it is proved that
\begin{thm}
If $f_k$ is an {\em RotS} function of degree $k$, then the nonlinearity
satisfies $N_{f_k}\geq 2^{n-k}.$
\end{thm}

Now, we evaluate the nonlinearity of $f_2$ for $n$ even.
\begin{lem}
\label{t_2k}
For $n\geq 3$, let
$t_n=x_1 x_2+x_2x_3+\cdots+x_{n-2}x_{n-1}+x_{n-1} x_{n}$.
 Then $t_{2k}$ is a bent function.
\end{lem}
\begin{proof}
We have
\[
\begin{split}
t_{2k}=
&x_2(x_1+x_3)+x_4(x_3+x_5)+\cdots +\\
&x_{2k-2}(x_{2k-3}+x_{2k-1})+x_{2k}x_{2k-1}.
\end{split}
\]
By taking the transformation
\[
X_{2i}=x_{2i}\ \text{and}\ X_{2i-1}=x_{2i-1}+x_{2i+1},
X_{2k-1}=x_{2k-1}, i=1,2,\ldots,k-1,
\]
we see that $t_{2k}$ is affinely equivalent to a bent function in
the Maiorana-McFarland class (see \cite{D1}), therefore it is also bent.
\end{proof}

We say (see  \cite{CLK}) that $g\in W_{2k+1}$ is {\em semi-bent},
if there is a bent
function $g_0\in W_{2k}$ with
\[
g=g_0||g_1,
\]
 where $g_1(x)=g_0(Ax\op a)\op 1$, $A$ is a nonsingular $2k$ by $2k$
  matrix and $a$ is any vector in $W_{2k}$.

In \cite{CLK}, the authors prove the following results (see Theorem 18,
 Corollary 21 and
Theorem 16), which will be used in this paper.
\begin{lem}
\label{clk1}
Any semi-bent function $g\in W_{2k+1}$ is balanced, $N_g=2^{2k}-2^k$,
 for any
$w^*\in W_{2k+1}$, the
correlation value between $g$ and the linear function
$l_{w^*}(x)=w^*\cdot x$ is $0$
or $\pm 2^{-k}$, and
\[
\#\{w^*\in W_{2k+1}|\, c(g,l_{w^*})=0\}=2^{2k}=\#\{w^*\in W_{2k+1}|\,
c(g,l_{w^*}=\pm 2^{-k}\}.
\]
\end{lem}

\begin{lem}
\label{clk2}
Let $g\in W_{2k+1}$ be a semi-bent function with $A=I$ and
 $a=(1,1,\cdots,1)$. Then $g$
satisfies $PC(2k)$.
\end{lem}

\begin{lem}
\label{clk3}
If $g$ is the concatenation  $g_0||g_1$, $w^*=(w,w_{2n+1})\in W_{n+1}$,
then
\[
\hF_{\hat g}(w^*)=\hF_{\hat g_0}(w)+(-1)^{w_{n+1}}\hF_{\hat g_1}(w).
\]
\end{lem}
The following result belongs to Preneel \cite{P}. We define
$l_b(x)=b\cdot x$.
\begin{lem}
\label{preneel}
For $h$ on $W_n$, $a,b\in W_n, c\in W_1$ and a $2k\times 2k$
nonsingular matrix $A$,
define $g$ by $g(x)=h(Ax\op a)\op l_b(x)\op c$. Then,
\[
\hF_{\hat g}(w)=(-1)^c(-1)^{(A^{-1}a,w\op b)}\hF_{\hat h}
((A^{-1})^t(w\op b)).
\]
\end{lem}

It is not very difficult to observe (see also \cite{PQ}) that any
2-degree rotation-symmetric function in $n$ variables
is affinely equivalent to
$f_2^n=f_2=x_1x_2\op x_2x_3\op \cdots \op x_{n-1}x_n\op x_nx_1.$
We show below that  $f_2^{2k}$  is not bent.
To do that we display an algorithm to evaluate $f_2^{2k}$ fast.
For that we need the following lemma, which can be proved by
considering the truth table.
\begin{lem}
Each monomial of degree 2 can be written in the form (\ref{f-blocks}) as
\begin{eqnarray}
x_i x_j &=& \left(D_{2^{n-i-2}}\left(D_{2^{n-j-2}}\bar
D_{2^{n-j-2}}\right)_{2^{j-i-1}}\right)_{2^{i-1}},\nonumber\\
& &{\rm if}\ 1\leq i<j\leq n-2,\nonumber\\
x_i x_{n-1} &=& \left(D_{2^{n-i-2}} A_{2^{n-i-2}}\right)_{2^{i-1}},
\label{2monom}\\
x_ix_n &=& \left(D_{2^{n-i-2}} B_{2^{n-i-2}}\right)_{2^{i-1}},\nonumber\\
x_{n-1}x_n &=& V_{2^{n-2}}.\nonumber
\end{eqnarray}
\end{lem}

Using (\refth{2monom}) we see that
\begin{eqnarray}
 f_2 &=& x_1x_2\op x_2x_3\op \cdots x_{n-1}x_n\op x_nx_1=\nonumber\\
 & & D_{2^{n-3}}\left(D_{2^{n-4}}\bar D_{2^{n-4}}\right)\op
 \left( D_{2^{n-4}}\left(D_{2^{n-5}}\bar D_{2^{n-5}}\right)
 \right)_2\op\\
& & \left( D_2\left(D\bar D\right) \right)_{2^{n-4}}\op
(DA)_{2^{n-3}}\op (DB)_{2^{n-3}}\op V_{2^{n-2}}=\nonumber\\
& & g\op (DC)_{2^{n-3}}\op V_{2^{n-2}}=g\op (V\bar U)_{2^{n-3}}\op D_{2^{n-3}}B_{2^{n-3}},\nonumber
\end{eqnarray}
where $g$ is the sum of the first $n-3$ strings of length $2^n$.

For a string $u$ of length $2^s, s\geq 4$, we denote by $\tilde u$,
the string obtained by complementing the second half, that is the
last $2^{s-1}$ bits of $u$.
It is not difficult to observe that the following algorithm will
output $f_2=G_1||G_2||G_3.$

{\bf Algorithm f2}.\\
\em
{\bf step $3$}:    $g_1^3\gets V Y,g_2^3\gets X\bar U$\\
{\bf step $s$}:  $g_i^s\gets g_i^{s-1}||\tilde g_i^{s-1},i=1,2$\\
  {\bf output}: $G_1\gets g_1^{n-4},G_2\gets g_2^{n-5}$,
        $G_3=\bar G_4$, where $G_4={\tilde {G_2}}$, and write $f_2=G_1||G_2||G_3$

\rm
For instance, the first three steps of the algorithm will produce
\begin{eqnarray*}
&&G_1\gets ((V Y)(V\bar Y))((V Y)(\bar V Y))\\
&&G_2||G_3\gets (X\bar U X U)(\bar X U X U).
\end{eqnarray*}

\begin{thm}
\label{f2}
If $f_2$ is defined on $W_n$, with $n=2k$, then
it is not bent.
Moreover, the nonlinearity
is
\[
N_{f_2}=2^{2k-1}-2^k,
\]
and the truth table of $f_2$ can be displayed using only $2^{n-3}-2$
operations (additions and multiplications).
\end{thm}

\begin{proof}
Using the above algorithm, we deduce that
the {\em RotS} function on $W_n$ of degree 2 can be be evaluated in $n-2$
steps, which requires
\[
(1+2^1+\cdots + 2^{n-5})+(1+2^1+\cdots+2^{n-6})+2^{n-5}=2^{n-3}-2
\]
operations, since at each step $s$ we complement $2^{s-2}$ bits.

First, we take an example, say $f_2^5=VYV\bar Y X \bar U\bar X\bar U=t_5+x_1x_5$
on $W_5$. We see that $f_2^5=t_4({\bf x})||t_4({\bf x}\op {\bf 1})\op 1$, therefore it is
semi-bent.

It is very easy to see that
\[
f_2^{2k+1}=t_{2k}({\bf x}_{2k})||(t_{2k}({\bf x}_{2k})+x_1+x_{2k}).
\]
But
\[
\begin{split}
\overline{t_{2k}({\bf x}_{2k})+x_1+x_{2k}}=
& \sum_{i=1}^{2k-1} x_ix_{i+1}+x_1+x_{2k}+1=\\
=&\sum_{i=1}^{2k-1}(x_i+1)(x_{i+1}+1)=t_{2k}({\bf x}_{2k}\op {\bf 1}),
\end{split}
\]
therefore $f_2^{2k+1}=t_{2k}({\bf x}_{2k})||(t_{2k}({\bf x}_{2k}\op {\bf 1})\op 1)$
is semi-bent. By Lemma \refth{clk2}, $f_2^{2k+1}$ satisfies the propagation criterion for all weights
$1\leq w\leq 2k$.

Similarly,
\[
f_2^{2k}=t_{2k-1}({\bf x}_{2k-1})||(t_{2k-1}({\bf x}_{2k-1})+x_1+x_{2k-1}).
\]

Now, we shall use Lemma \refth{clk3} to compute the nonlinearity of $f_2^{2k}$.
First, we observe that
\[
t_{2k+1}=t_{2k}({\bf x}_{2k})||(t_{2k}({\bf x}_{2k})+x_{2k}).
\]
Take $A=I$ and $a=(1,0,1,0,\ldots,1,0)$. We see that
\begin{eqnarray*}
t_{2k}({\bf x}_{2k})+x_{2k}&=&t_{2k}(x_1+1,x_2,x_3+1,\ldots,x_{2k-1},x_{2k})\\
&=&(x_1+1)x_2+x_2(x_3+1)+\cdots+(x_{2k-1}+1)x_{2k}.
\end{eqnarray*}
We denote the last expression by $r({\bf x})$.
Using  Lemma \refth{clk3}  we compute the Walsh-Hadamard transform
\[
\hF_{\hat r}({\bf w}_{2k})=(-1)^{({\bf w},{\bf a})}\hF_{\hat t_{2k}}({\bf w})=\pm 2^k,
\]
since by Lemma \refth{t_2k}, $t_{2k}$ is bent.

For simplicity we set $t({\bf x})=t_{2k+1}({\bf x}_{2k+1})$
and $w^*=(w,w_{2k+1})$. Thus,
\beq
\hF_{\hat t}(w^*)=\hF_{\hat t_{2n}}({w})+(-1)^{w_{2k+1}}\hF_{\hat r}(w)=0\ \text{or}\ \pm 2^{k+1},
\end{equation}
since $r$ and $t_{2k}$ are bent. Therefore,
\[
N_{t_{2k+1}}=2^{2k}-\frac{1}{2}|\hF_{\hat t_{2k+1}}(w^*)|=2^{2k}-2^k.
\]
By Lemma \refth{clk1}, $\hF_{\hat t_{2k-1}}({\bf x}_{2k-1})=0$ or $\pm 2^k$. Let
$v({\bf x})=t_{2k-1}({\bf x})+x_1+x_{2k-1}$. By Lemma \refth{clk3},
\[
\hF_{\hat v}({\bf x})=(-1)^{({\bf x}\op (1,0,\ldots,0,1),{\bf 0})} \hF_{\hat t_{2k-1}}\left({\bf x}\op (1,0,\ldots,0,1)\right)=0\ \text{or}\ \pm 2^k.
\]
Thus, by the same Lemma \refth{clk3},
\[
\hF_{\hat f_2^{2k}}({\bf x}_{2k})=\hF_{\hat t_{2k-1}}({\bf x}_{2k-1})+(-1)^{x_{2k}} \hF_{\hat v}({\bf x})=0\ \text{or}\ \pm 2^{k+1},
\]
which implies $N_{f_2^{2k}}=2^{2k-1}-2^k$.
Therefore $f_2$ is not bent (any bent function in $2k$ variables has
nonlinearity $2^{2k-1} - 2^{k-1}$  [2, Th. 13, p. 111]) and the
theorem is proved.
\end{proof}

\rm
\begin{rem}
We remark that, using the normal form of the function, the truth table
of $f$ is found using $\frac{3n-1}{2} 2^n$ operations (see
\cite{PQ} for a detailed discussion). Using the previous theorem
we can display the truth table using only $2^{n-3}-2$ operations, which is
 a significant improvement.
\end{rem}

Now, we will evaluate the weights of $f_2$ for any dimension $n$.
We prove
\begin{thm}
The weights of $f_2$ are given by
\begin{equation}
wt(f_2^n)=2^{n-1}-2^{\frac{n}{2}-1}\left(1+(-1)^n\right).
\end{equation}
\end{thm}

\begin{proof}
We recall that $f_2=g_1^{n-1} g_2^{n-2} g_3^{n-2}$. We show that
for any $s$,
\begin{equation}
\label{weights_f2}
wt(g_i^s)=2 wt(g_i^{s-2})+2^{s-2}, i=1,2,3.
\end{equation}

Since $g_i^s=g_i^{s-1}\tilde g_i^{s-1}=g_i^{s-1}g_i^{s-2}\bar {\tilde g}_i^{s-2}$,
\begin{eqnarray}
wt(g_i^s)&=&wt(g_i^{s-1})+wt(g_i^{s-2})+wt(\tilde g_i^{s-2})\nonumber\\
&=& wt(g_i^{s-2})+wt(\tilde g_i^{s-2})+wt(g_i^{s-1})+2^{s-2}-wt(\tilde g_i^{s-2})\label{gi}\\
&=& 2\, wt(g_i^{s-2}) +2^{s-2}, i=1,2.\nonumber
\end{eqnarray}

Now, from $g_3^s=\bar g_2^{s-1}\tilde g_2^{s-1}$, we get
\begin{eqnarray}
wt(g_3^s)&=&2^{s-1}-wt(g_2^{s-1})+wt(\tilde g_2^{s-1})\nonumber\\
&=&2^{s-1}-wt(g_2^{s-1})+2wt(g_2^{s-2})-wt(g_2^{s-1})+2^{s-2} \label{g3_1}\\
&=& 2wt(g_2^{s-2})- 2wt(g_2^{s-1})+ 2^{s-1}+2^{s-2} \nonumber\\
&=&wt(g_2^{s})-2wt(g_2^{s-1})+2^{s-1}.\nonumber
\nonumber
\end{eqnarray}
The above equation, for $s-1$, produces
\begin{equation}
\label{g3_2}
wt(g_3^{s-1})=wt(g_2^{s-1})-2wt(g_2^{s-2})+2^{s-2}.
\end{equation}
Now, we add \refeq{g3_1} plus twice \refeq{g3_2}, and we get
\[
wt(g_3^s)+2 wt(g_3^{s-1})=wt(g_2^{s})-4wt(g_2^{s-2})+2^s.
\]
But $wt(g_2^{s})=2 wt(g_2^{s-2})+2^{s-2}$. By adding the two previous
equations we get
\begin{equation}
\label{g3_3}
wt(g_2^{s-2})=2^{s-1}+2^{s-3}-wt(g_3^{s-1})-\frac{wt(g_3^{s})}{2}
\end{equation}
Replacing \refeq{g3_3} into \refeq{g3_2} we obtain
\[
wt(g_3^{s+2})=2wt(g_3^s)+2^s.
\]
This together with \refeq{gi} will give the following recurrence for the weights of $f_2$.
\begin{equation}
wt(f_2^n)=2 wt(f_2^{n-2})+2^{n-2}.
\end{equation}

A generating function for the above recurrence is
\begin{equation}
- {\displaystyle \frac {32\,{\displaystyle \frac {z^{7}}{1 - 2\,z}}  +
16\,z^{5} + 24\,z^{6}}{ - 1 + 2\,z^{2}}}.
\end{equation}

We can linearize the recurrence by using the transformation
\[
y_n=wt(f_2^n)-2^{n-1},
\]
thus obtaining the recurrence
\[
y_{n}=2y_{n-2}.
\]
Using the above simple recurrence with $wt(f_2^5)=16$ and
$wt(f_2^6)=24$ we get
a closed formula for the weights of $f_2$ in dimension $n$,
namely
\[
 2^{n-1}-2^{\frac{n}{2}-1}\left(1+(-1)^n\right),
\]
and the theorem is proved.
\end{proof}


\section{The third degree rotation-symmetric function}


\rm

As in the case of second degree {\em RotS} functions, it is easy to observe
that any {\em RotS} function of degree 3
in $n$ variables, $f_3^n=f_3$ is affinely equivalent to
\begin{equation}
f_3 = x_1 x_2 x_3+ x_2 x_3 x_4+\cdots+x_n x_1 x_2.
\end{equation}
Now, using a computer program we have determined the
nonlinearity of $f_3$ on $W_n,n\geq 9$ which turns out to be the same as its weight.
Thus,
\[
\begin{tabular}{|c|c|c|c|c|c|c|c|}\hline
$n$ & $3$ & $4$ & $5$ & $6$ & $7$ & $8$ & $9$\\\hline
 $N_{f_3^n}$ & $1$ &  $4$ &  $6$ &  $18$ &  $36$ &  $80$ &  $172$\\\hline
\end{tabular}
\]
We shall assume that $n\geq 10$.
The following lemma will be used.

\begin{lem}
\label{s_monomial}
The truth table of any monomial $x_{i_1} \cdots x_{i_s}$ of degree $s$ is
\begin{eqnarray}
&& \left(D_{2^{n-i_1-2}}\cdots\left(D_{2^{n-i_s-2}}
\bar D_{2^{n-i_s-2}}\right)_{2^{i_s-i_{s-1}-1}} \right)_{2^{i_1-1}},\nonumber\\
&&\ { if}\ 1\leq i_1<\cdots<i_s\leq n-2,\nonumber\\
&& \left(D_{2^{n-i_1-2}}\cdots\left(D_{2^{n-i_{s-1}-2}}
M_{2^{n-i_{s-1}-2}}\right)_{2^{i_{s-1}-i_{s-2}-1}} \right)_{2^{i_1-1}},\label{s_monom}\\
&&\ {where}\ M=A\ { or}\ B\ { if}\ i_s=n-1,\ {respectively}\ i_s=n,\nonumber\\
&& \left(D_{2^{n-i_1-2}}\cdots\left(D_{2^{n-i_{s-2}-2}}
V_{2^{n-i_{s-2}-2}}\right)_{2^{i_{s-2}-i_{s-3}-1}}\right)_{2^{i_1-1}},\nonumber\\
&&\ { if}\ i_{s-1}=n-1 \ {and}\ i_s=n.\nonumber
\end{eqnarray}
\end{lem}
\begin{proof}
Straightforward using the truth table.
\end{proof}

Using the above lemma we write
\begin{equation}
x_i x_{i+1} x_{i+2}=\left(D_{2^{n-i-2}}\left(D_{2^{n-i-3}}
\left(D_{2^{n-i-4}}\bar
D_{2^{n-i-3}}\right)\right)\right)_{2^{i-1}},
\end{equation}
if $i\leq n-4$, and
\begin{eqnarray}
x_{n-3}x_{n-2} x_{n-1} &=& (D_3 A)_{2^{n-4}}\nonumber\\
x_{n-2}x_{n-1} x_{n} &=& (D V)_{2^{n-3}}\\
x_{n-1}x_{n} x_{1} &=& D_{2^{n-3}} V_{2^{n-3}}\nonumber\\
x_{n}x_{2} x_{1} &=& D_{2^{n-3}+2^{n-4}} B_{2^{n-4}}\nonumber
\end{eqnarray}
Therefore,
\begin{eqnarray*}
f_3
&=& \sum_{i=1}^{n-4}\left(D_{2^{n-i-2}+2^{n-i-3}+2^{n-i-4}}\bar
D_{2^{n-i-4}}\right)_{2^{i-1}}\oplus  (D_3 A)_{2^{n-4}}\op\\
&& (D V)_{2^{n-3}}\op  D_{2^{n-3}+2^{n-4}} B_{2^{n-4}}\op
D_{2^{n-3}} V_{2^{n-3}}=\\
&=& \sum_{i=1}^{n-4}\left(D_{2^{n-i-2}+2^{n-i-3}+2^{n-i-4}}\bar
D_{2^{n-i-4}}\right)_{2^{i-1}}\oplus\\
&& (D V D Y)_{2^{n-4}}\op D_{2^{n-3}} V_{2^{n-4}} X_{2^{n-4}}=\\
&=&  \sum_{i=1}^{n-4}\left(D_{2^{n-i-2}+2^{n-i-3}+2^{n-i-4}}\bar
D_{2^{n-i-4}}\right)_{2^{i-1}}\op\\
&&
\left( D_3 C\right)_{2^{n-5}}
\left( V_3\bar U\right)_{2^{n-6}}
\left( X_3 Y\right)_{2^{n-6}}=H_1||H_2||H_3||H_4,
\end{eqnarray*}
where $H_1$ (on $W_{n-1}$), $H_2$ (on $W_{n-2}$), $H_3,H_4$
(on $W_{n-3}$) are
defined by the following algorithm ($\hat u$, on $W_j$, is the string obtained
from $u$ by complementing its last $2^{j-2}$ bits):

\newpage

{\bf Algorithm f3.}\\
\em
{\bf step $4$}:
$h_1^4\gets DVDY$  $h_2^4\gets VDVA$, $h_3^4\gets XBXC$\\
{\bf step $s$}:
$h_i^s\gets h_i^{s-1}||\hat h_i^{s-1}$\\
{\bf output}: $H_1\gets h_1^{n-1},H_2\gets h_2^{n-2},H_3\gets
h_3^{n-3}$, $H_4$ is the string obtained from $\hat H_3$ by
complementing its first half, that is $H_4=\bar H_5$, where $H_5=\tilde H_6, H_6=\hat H_3.$
Write $f_3=H_1||H_2||H_3||H_4.$
\rm

As in the case of the 2-degree {\em RotS} function we see that we need
\[
\begin{split}
& 2^2(1+2+\cdots 2^{n-5})+2^2(1+2+\cdots 2^{n-6})+\\
&2^2(1+2+\cdots 2^{n-7})+2^{n-4}+2^{n-5}= 3\cdot 2^2(2^{n-6}-1)+\\
& 2^{n-3}+2^{n-4}=2^{n-2}+2^{n-4}+2^{n-5}-3\cdot 2^2
\end{split}
\]
operations to display the truth table of $f_3^n$.

We shall  evaluate the weight of $f_3^s$ for any $s$. To do this we will compute
the weights of each component of $f_3^s$. We observe that
\[
\begin{split}
 h_i^s&=h_i^{s-1}h_i^{s-2} h_i^{s-3}\bar h_i^{s-4} \hat {\bar h}_i^{s-4}\ \text{and}\\
 \hat h_i^s&=h_i^{s-1}h_i^{s-2}\bar h_i^{s-3}h_i^{s-4} \hat {h}_i^{s-4}, i=1,2,3.
\end{split}
\]
Therefore, denoting by $w_i^s$ the weight of $h_i^s$, and by  $\hat w_i^s$
the weight of $\hat h_i^s,i=1,2,3$,
 we arrive at the following identities:
\begin{eqnarray}
\hat w_i^s&=&2w_i^{s-1}+2w_i^{s-2}-w_s+2^{s-2},\label{hat_weight}\\
w_i^s&=&w_i^{s-1}+\hat w_i^{s-1}.\label{weight}
\end{eqnarray}

Using Mathematica\footnote{A trademark of {\it Wolfram Research}} we obtained the
following results on the weights of $f_3^n$ and of each of the four
components on dimensions less than 12.
\beq
\label{table_weights}
\begin{tabular}{|l|c|c|c|c|c|}\hline
{\em n}& $wt(f_3^n)$  & $wt(h_1^{n-1})$   &  $wt(h_2^{n-2})$  &
$wt(h_3^{n-3})$  &  $wt(h_4^{n-3})$                  \\\hline
$3$ & $1$ & & & &\\\hline
$4$ & $4$ & & & &\\\hline
$5$ & $6$ & $2$  & & & \\\hline
 $6$& $18$ & $6$  & $4$  & &\\\hline
 $7$& $36$ & $14 $ & $8$ & $6$  &  $8$      \\\hline
 $8$& $80$ & $32$  & $18$ & $12$  &  $18$      \\\hline
$9$& $172$ & $72$  & $40$ & $26$  &  $34$    \\\hline
$10$& $360$ & $156$  & $84$ & $52$  &  $68$      \\\hline
$11$& $760$ & $336$  & $180$ & $108$  &  $136$      \\\hline
$12$& $1576$ & $712$  & $376$ & $220$  &  $268$      \\\hline
\end{tabular}
\end{equation}
We have
\[
wt(f_3^n)=wt(h_1^{n-1})+wt(h_2^{n-2})+wt(h_3^{n-3})+wt(h_4^{n-3}).
\]

We show by induction that
\beq
\label{rec1}
wt(h_i^s)=2\left(wt(h_i^{s-2})+wt(h_i^{s-3})\right)+2^{s-4}, i=1,2,3,4.
\end{equation}
From the table \refeq{table_weights} we have the truth of the claim for the first few cases.
Assume \refeq{rec1} true for $s-1$ and we prove it for $s$.
From \refeq{hat_weight} and \refeq{weight} and by using the induction step we get
\[
wt(h_i^s)=wt(h_i^{s-1})+wt(\hat h_i^{s-1})=2 \left(wt(h_i^{s-2})+wt(h_i^{s-3})  \right) +2^{s-3}.
\]
Similarly for $h_4^s$. Adding these relations we get
\beq
\label{rec_f3}
wt(f_3^s)=2\left( wt(f_3^{s-2})+wt(f_3^{s-3}) \right)+2^{s-3}.
\end{equation}
Remark that this equation is true for any $s\geq 6$.

Using the table \refeq{table_weights}, the recurrence \refeq{rec_f3} and
Maple\footnote{A trademark of {\it Waterloo Maple}}, we  get
\begin{thm}
The  generating function for the weight of $f_3$, is
\begin{equation}
\label{gen_f3}
 - {\displaystyle \frac {8\,{\displaystyle \frac {z^{
6}}{1 - 2\,z}}  + z^{3} + 4\,z^{4} + 4\,z^{5}}{ - 1 + 2\,z^{2} +
2\,z^{3}}}.
\end{equation}
\end{thm}

The series expansion of the above generating function is
\begin{eqnarray*}
&&z^{3} + 4\,z^{4} + 6\,z^{5} + 18\,z^{6} + 36\,z^{7} + 80\,z^{8}+\\
 &+& 172\,z^{9} + 360\,z^{10} + 760\,z^{11} + 1576\,z^{12} + {\rm O}(z^{13}),
\end{eqnarray*}
obtaining once again the weights of $f_3^n$, for any dimension.

Based on our numerical examples, we give the following conjecture.
\begin{conj}
The nonlinearity of $f_3^n$ is the same as its weight.
\end{conj}

{\noindent\sms Thomas W. Cusick:} {\em
State University of New York at Buffalo, Department of Mathematics,
Buffalo, NY 14260-2900, e-mail: cusick@math.buffalo.edu}\\
{\sms Pantelimon St\u anic\u a}: {\em
Auburn University Montgomery, Department of Ma\-the\-ma\-tics,
Montgomery, AL 36124-4023, e-mail: stanpan@strudel.aum.edu}

\end{document}